\documentclass[11pt]{article}

\usepackage{amsmath,amssymb,amsthm}
\usepackage{a4wide}
\usepackage[cp866]{inputenc}
\usepackage{url}
\usepackage{srcltx}

\allowdisplaybreaks

\def\nfrac#1#2{\mbox{\small$\dfrac{#1}{#2}$}}
\def\nbinom#1#2{\mbox{\small$\dbinom{#1}{#2}$}}

\newtheorem{theorem}{Theorem}

\begin{document}

\title{Fractional Matchings in Hypergraphs\thanks{Supported by
FAPESP, project nos.~2012/13341-8 and~2013/07699-0, and NUMEC/USP, project
MaCLinC/USP.}}


\author{Vladimir Blinovsky}

\date{\small Instituto de Matem\'atica e Estat\'{\i}stica, Universidade de S\~ao
Paulo,\\ Rua do Mat\~ao 1010, 05508-090, S\~ao Paulo, Brazil\\ Kharkevich Institute
for Information Transmission Problems,\\ Russian Academy of Sciences, B. Karetnyi 19,
Moscow, Russia\\ {\tt vblinovs@yandex.ru}}

\maketitle\bigskip

\begin{abstract}
We find an exact formula for the minimum number of edges in a hypergraph which
guarantees a fractional matching of cardinality $s$ in the case where $sn$ is an
integer.
\end{abstract}

\section{Introduction}

Let $\mathcal{H}=([n],E)$ be a $k$-uniform hypergraph with vertex set $[n]$ and a set
of edges $E\subset\nbinom{[n]}{k}$. Taking into account a natural bijection between
the set of binary $n$-tuples and the set $2^{[n]}$, we identify them in what follows.

A fractional matching of a hypergraph of cardinality $s\in [0,1]$ is a set of
nonnegative real numbers $\{\alpha_{e},\: e\in E\}$ such that $\sum\limits_{e\in
E}\alpha_e=s$ and the $n$-tuple $\bar{a}=(a_1,\ldots,a_n)=\sum\limits_{e\in
E}e\alpha_e$ has coordinates satisfying the inequalities $0\le a_j\le k/n$.

If $s\le{k}/{n}$, then the only hypergraph that has no fractional matching of
cardinality $s$ is the hypergraph without edges.

A fractional matching in the case $s=1$ is called a perfect fractional matching. This
case was considered in the paper~\cite{1}, where the following result was proved.

\begin{theorem}
The minimum number $M+1$ of edges in a hypergraph guaranteeing a perfect fractional
matching satisfies the equality
$$
M+1=\max_{n-1\ge a\ge 1}\sum_{i>ka/n}\binom{a}{i}\binom{n-a}{k-i} +1.
$$
\end{theorem}

This theorem was preceded by a conjecture formulated by Ahlswede and Khachatrian
in~\cite{2}.

In the present paper we find a formula for the minimum number of edges in a
hypergraph which has a fractional matching of cardinality $s$ in the case where $sn$
is an integer. As follows from the above, we may assume that $1>s>k/n$; we also
assume that $sn$ is an integer. We prove the following statement.

\begin{theorem}\label{th2}
The maximum number of edges $M(s,n,k)$ in a hypergraph which has no fractional
matching of cardinality $s$ satisfies the equality
\begin{equation}\label{e3}
M(s,n,k)=\max_{1\le c\le ns-1}\sum_{i>kc/ns}\binom{c}{i}\binom{n-c}{k-i}.
\end{equation}
\end{theorem}

In~\cite{3} the reader can find asymptotics of the function $M(s,n,k)$ as
$n\to\infty$ for several particular choices of $k$ and $s$.

\section{Proof of Theorem~\ref{th2}}

Let $\beta(n,k)\subset\mathbb{R}^n$ be a hypersimplex, i.e., a convex polytope with
the set of vertices $\smash[t]{\nbinom{[n]}{k}}$. In~fact, we are interested in a
transformed hypersimplex $s\beta(n,k)$ where each vector from $\beta(n,k)$ is
multiplied by $s$. Below we consider only such transformed hypersimplex. If a
hypergraph $\mathcal{H}=([n],E)$ has a fractional matching of cardinality $s$, then
the convex hull $X(E)$ of the vertices of this hypergraph in $s\beta(n,k)$ contains a
point $\bar{a}=(a_1,\ldots,a_n)$ whose coordinates belong to the interval $[0,k/n]$.
The set $A$ of all such points is also convex. This means that $\mathcal{H}$ has no
fractional matching of cardinality $s$ if and only if
\begin{equation}\label{e6}
X(E)\cap A=\varnothing.
\end{equation}
Thus, the original problem is reduced to the problem of finding the maximum
cardinality of a set $sE\in c\nbinom{[n]}{k}$ such that condition~\eqref{e6} holds.
If~\eqref{e6} holds, then there exists a hyperplane~$L$ such that $X(E)$ and $A$
belong to different half-spaces into which $L$ divides $\mathbb{R}^n$. Without loss
of generality we may assume that $L(\bar{0})=0$, $L(e)>0$ for $e\in E$, and
\begin{equation}\label{e8}
L(\bar{a})\le 0
\end{equation}
for $\bar{a}\in A$.

Condition~\eqref{e8} is equivalent to the condition that~\eqref{e8} is true for all
vertices $\bar{a}$ of the convex polygon $A\cap S$, where
$$
S=\{\bar{x}\in\mathbb{R}^n:\: (\bar{x},\bar{1})=sk,\: x_j\ge 0\}.
$$
Let $L=\{\bar{x}\in\mathbb{R}^n:\: (\bar{x},\bar{\omega})=0\}$ be the hyperplane
defined above. We assume that the coordinates $(\omega_1,\ldots,\omega_n)$ of the
hyperplane $L$ are not increasing: $\omega_1\ge\ldots\ge\omega_n$. Vertices of $A\cap
S$ are those which have $sn$ coordinates equal to $k/n$ and all other zero. We may
assume that the vector $\bar{a}=(k/n,\ldots,k/n,0,\ldots 0)$ belongs to the
hyperplane: $L(\bar{a})=0$. Thus, we have $\sum\limits_{i=1}^{sn}\omega=0$. We may
also assume that $w_i=w_{sn}$ for $i>sn$. The space of such $n$-tuples has a natural
basis $z_j=(sn-j,\ldots,sn-j,-j,\ldots,-j)$ for $j\in[sn-1]$, where the $j$th
vector~$z_j$ has $j$ coordinates $sn-j$. Any vector in this space is a linear
combination of the basis vectors with nonnegative coordinates. Let
$\bar{y}=\smash[b]{\sum\limits_{j=1}^{sn-1} z_j y_j}$ and $y_j\ge 0$. Then for
$\bar{x}=se\in s\beta(n,k)$ we have
$$
\begin{aligned}
(\bar{y},\bar{x})=s(\bar{y},e) &=s\sum_{i=1}^{n}e_i \sum_{j=1}^{sn-1}z_{ji}y_j\\
&=s\sum_{j=1}^{sn-1}y_j \sum_{i=1}^{n}e_i z_{ji}\\ &=s\sum_{j=1}^{sn-1}y_j
\,\Biggl(sn\sum_{i=1}^jx_i -jk\Biggr)\\ &=s\,\Biggl(sn \sum_{j=1}^{sn-1}y_j
\sum_{i=1}^jx_i -k\sum_{j=1}^{sn-1}jy_j\Biggr).
\end{aligned}
$$
Dividing the last expression in this chain of equalities by
$s\sum\limits_{j=1}^{sn-1}jy_j$ and imposing the condition of positiveness
of the scalar product, we obtain the inequality
$$
\sum_{j=1}^{sn-1}\frac{y_j}{\sum\limits_{i=1}^{sn-1}iy_i}\sum_{\ell=1}^j x_\ell
>\frac{k}{ns}.
$$
This is equivalent to the inequality
\begin{equation}\label{e0}
\sum_{j=1}^{sn-1}\alpha_j x_j >\frac{k}{sn}
\end{equation}
for some $\alpha_j\ge 0$ such that
$\smash{\sum\limits_{i=1}^{sn-1}\alpha}=1$.\rule{0pt}{12pt} Hence, to complete the
proof of Theorem~\ref{th2}, we have to show that the maximum (over the choices of
$\alpha$) number of solutions in~$\nbinom{[n]}{k}$ of the inequality~\eqref{e0} is
$M(s,n,k)$. To prove this, we use the technique from~\cite{1}. Consider the function
$$
f(\{\alpha_1,\ldots,\alpha_{sn-1}\})=\frac{1}{\sqrt{2\pi}}\sum_{x\in\binom{[n]}{k}}
\int_{-\infty}^{\frac{\sum\limits_{j=1}^{sn-1}\alpha_j x_j
-\frac{k}{sn}}{\sigma}}e^{-\frac{z^2}{2}}dz
$$
Define
$$
N(\alpha_1,\ldots,\alpha_{sn-1})=\Biggl|x\in\binom{[n]}{k}:\:
\sum_{j=1}^{sn-1}\alpha_j x_j >\frac{k}{sn}\Biggr|.
$$
Then we have
$$
|N(\{\alpha_j\})-f(\{\alpha_j\})|<\epsilon (\sigma),\quad
\epsilon(\sigma)\xrightarrow{\sigma\to 0\,}0,
$$
uniformly over $\{\alpha_j\}$ such that
\begin{equation}\label{e44}
\Biggl|\sum_{j=1}^{sn-1}\alpha_j x_j -\frac{k}{sn}\Biggr|>\delta,\quad \text{for
all}\ x\in\binom{[n]}{k}.
\end{equation}
For the extremal $\alpha$ with $N(\alpha)=M(s,k,n)$, it is easy to see that~$\alpha$
satisfies condition~\eqref{e44} for some $\delta>0$, because otherwise, if
$\smash[b]{\sum\limits_{j=1}^{sn-1}\tilde{\alpha}_j x^0_j}=\smash[b]{\nfrac{k}{sn}}$
for some $x^0\in\smash[b]{\nbinom{[n]}{k}}$, then for $\alpha'$ sufficiently close
to $\tilde{\alpha}$ the conditions
$$
\alpha_j'\ge 0,\qquad \sum_{j=1}^{sn-1} \alpha_j'=1,\qquad \sum_{j=1}^{sn-1}\alpha'_j
x_j >\frac{k}{sn}\rule{0pt}{22pt}
$$
are not violated. Hence, when we are interested in extremal $\alpha$, we may assume
that~\eqref{e44} is satisfied.

We will assume without loss of generality that $\alpha_1\ge\ldots\ge\alpha_{sn-1}$.
Since we have the restrictions $\alpha_j\ge 0$, we should look for the extremum among
$\alpha$ such that
$$
\alpha_{a+1}=\ldots=\alpha_{sn-1}=0,\quad a\in [sn-1]
$$
(the case $a=sn-1$ means that we are not imposing any zero condition on $\alpha$).
Assume that this condition is valid for some $a$. Then, because
$\alpha_a=1-\smash[b]{\sum\limits_{j=1}^{a-1}\alpha_j}$, we
have\rule[-4pt]{0pt}{10pt}
\begin{equation}\label{e45}
f'_{\alpha_j}=\frac{1}{\sqrt{2\pi}\sigma}\sum_{x\in\binom{[n]}{k}:\: j\in x,\:
a\notin x}e^{-\frac{\bigl(\sum\limits_{j=1}^{a}\alpha_j x_j
-\frac{k}{sn}\bigl)^2}{2\sigma^2}}-\frac{1}{\sqrt{2\pi}\sigma}
\sum_{x\in\binom{[n]}{k}:\: j\notin x,\: a\in x}
e^{-\frac{\bigl(\sum\limits_{j=1}^{a}\alpha_j x_j -\frac{k}{sn}\bigr)^2}{2\sigma^2}}.
\end{equation}
In what follows, we assume that $a>4$. The cases $a\le 4$ are easy to
treat.

Note that if $\sum\limits_{j=1}^{a}x_j=\nfrac{ka}{sn}$ for some
$x\in\nbinom{[n]}{k}$, then $f(\alpha)$ does not achieve its (global) extremum on
$\alpha=\Bigl(\nfrac{1}{a},\ldots,\nfrac{1}{a},0,\ldots,0\Bigr)$ when $\sigma$ is
small. This can be shown using the same small perturbation arguments as above.

Now let us show that we may assume that these equalities can be valid together on
step functions $\beta_j=\beta_a$ for $j\in [a]$. Indeed, choose the parameter
$\sigma$ sufficiently small and then fix~it. Then, to satisfy equations~\eqref{e45},
we should assume that the equalities
\begin{eqnarray*}
\sum_{x\in\binom{[n]}{k}:\: j\in x,\: a\notin x}
e^{-\frac{\left((\beta,x)-\frac{k}{sn}\right)^2}{2\sigma^2}} &=&\nonumber
\sum_{x\in\binom{[n]}{k}:\: a\in x,\: j \notin x}
e^{-\frac{\left((\beta,x)-\frac{k}{sn}\right)^2}{2\sigma^2}}
\end{eqnarray*}
are valid. To satisfy these equalities, we should assume that the exponents in the
sums on the left- and right-hand sides are equal; i.e., for each given $j\in [a -1]$
\begin{equation}\label{er89}
\Bigl((\beta,x)-\frac{k}{sn}\Bigr)^2=\Bigl((\beta,y)-\frac{k}{sn}\Bigr)^2
\end{equation}
where $x\in\nbinom{[n]}{k}$, $j\in x$, $y\in\nbinom{[n]}{k}$, $a\in y$, and
$x\setminus j$ and $y\setminus a$ run over all sets of cardinality $k-1$ in
$[n-j-a]$. We rewrite equalities~\eqref{er89} as follows:
\begin{multline*}
\beta_j^2 +(\beta_{j_1}+\ldots+\beta_{j_{k-1}})^2 -2\frac{k}{sn} \beta_j
-2\frac{k}{sn} (\beta_{j_1}+\ldots+\beta_{j_{k-1}})
+\beta_j(\beta_{j_1}+\ldots+\beta_{j_{k-1}})\\ = \beta_a^2
+(\beta_{m_1}+\ldots+\beta_{m_{k-1}})^2 -2\frac{k}{sn} \beta_a-2\frac{k}{sn}
(\beta_{m_1}+\ldots+\beta_{m_{k-1}}) +\beta_a(\beta_{m_1}+\ldots+\beta_{m_{k-1}}).
\end{multline*}
Summing up both sides of these equality over all admissible choices of
$j_1,\ldots,j_{k-1}$ and $m_1,\ldots,m_{k-1}$ leads to the equality
\begin{equation}\label{er56}
\binom{n-2}{k-1}\left(\beta_j^2 -2\frac{k}{sn} \beta_j\right) -2\frac{k}{sn} R
+2\beta_jR =\binom{n-2}{k-1}\left(\beta_a^2 -2\frac{k}{sn} \beta_a\right)
-2\frac{k}{sn} R +2\beta_aR,
\end{equation}
where
$$
R= \sum_{x\in\binom{[n]\setminus\{j, a\}}{k-1}} (\beta,x)=\binom{n-3}{k-2}\sum_{m\ne
j,a}\beta_m =\binom{n-3}{k-2}(1-\beta_j-\beta_a).
$$
From~\eqref{er56} it follows that $\beta_j$ can take at most two values:
\begin{equation}\label{er91}
\begin{aligned}
\beta_j&=\beta_a,\\[-3pt] \beta_j+\beta_a&=\lambda \triangleq 2\frac{\nfrac{k}{sn}
-\nfrac{k-1}{n-2}}{1-2\nfrac{k-1}{n-2}}.
\end{aligned}
\end{equation}
Next we show how we can eliminate the possibility that $\beta_j$ takes the second
value. First assume that to each $x$ such that $|x\cap [a]|=p$ there corresponds some
$y$ such that $|y\cap [a]|=p$ for all $x\in\nbinom{[n]}{k}$ and $p$. For a given~$p$
we sum up the left- and right-hand sides of~\eqref{er89} over~$x$ and the
corresponding $y$ such that $|x\cap [a]|=p$. Then, similarly to the case of summation
over all $x$, we obtain two possibilities: either
$$
\beta_j=\beta_a
$$
or
\begin{equation}\label{et1}
\beta_j+\beta_a=2\frac{\nfrac{k}{sn} -\nfrac{p-1}{a-2}}{1-2\nfrac{p-1}{a-2}}.
\end{equation}
Since $p$ can be varied, it follows that the last equality for some $p$ contradicts
the second equality in~\eqref{er91}.

Now assume that for some $b$ we have
\begin{equation}\label{er77}
\beta_j=
\begin{cases}
\lambda-\beta_a,& j\le b,\\ \beta_a,& j\in [b+1,a].
\end{cases}
\end{equation}
Since $\sum\limits_j \beta_j=1$, we have the following condition on $\beta_a$ and
$\nfrac{k}{sn}$:
\begin{equation}\label{ek1}
b\lambda +(a -2b)\beta_a=1.
\end{equation}
Let $\beta_j=\lambda -\beta_a$. Assume also that to some $x$ such that $|x\cap
[a]|=p$ there corresponds some~$y$ such that $|y\cap [a]|=q$ for some $p\ne q$.
From~\eqref{er89} it follows that there are two possibilities: either
$$
(\beta,x)=(\beta,y)
$$
or
\begin{equation}\label{e34}
(\beta,x)+(\beta,y)=2\frac{k}{sn}.
\end{equation}
Each of these equalities impose some condition; the first equality, the condition
(for some integers $p_1$ and $p_2$)
$$
p_1\beta_a+p_2\lambda=0,
$$
which is either inconsistent with equality~\eqref{ek1} or together with
equality~\eqref{ek1} uniquely determines the value of $\nfrac{k}{sn}$.

On the other hand, equality~\eqref{e34} imposes the condition (for some integers
$p_3,p_4$)
\begin{equation}\label{ed1}
p_3 \beta_a+p_4 \lambda=2\frac{k}{sn}.
\end{equation}
It is possible that equality~\eqref{ek1} together with equality~\eqref{ed1} does not
determine the value of~$k/n$. In~this case there again can be two possibilities. The
first is that there exist $x$ such that $|x\cap [a]|=m$ (where $m$ can be equal to
either $p$ or $q$) and the corresponding~$y$ such that $|y\cap [a]|=v$ with $v\ne
p,q$.

The second possibility is that to each $x$ such that $|x\cap [a]|=m$ with $m\ne p,q$
there corresponds some $y$ such that $|y\cap [a]|=m$. In this second case we again
come to the case that leads to equalities~\eqref{et1} (because for $a\ge 5$ the
number of such $m\ne p,q$ is greater than~$1$).

If we have the first possibility, then there is an additional equation
\begin{equation}\label{el1}
q_3 \beta_a+q_4 \lambda=2\frac{k}{sn}
\end{equation}
which together with~\eqref{ek1} and~\eqref{ed1} is either inconsistent or determines
a unique value of $\smash[t]{\nfrac{k}{sn}}$.

We see that if $b>1$ and $\beta_j=\beta -\beta_a >\beta_a$ for $j\le b$, then $\beta$
can take values only in some discrete finite set. Varying the
value~$\smash[b]{\nfrac{k}{sn}}$ a little (considering instead of
$\smash[b]{\nfrac{k}{sn}}$ other numbers sufficiently close to
$\smash[b]{\nfrac{k}{sn}}$), we can achieve the situation where neither of values of
these functions coincides with the true value of $\nfrac{k}{sn}$. Again we note that
such small perturbation can always be done without violating relation~\eqref{e44}.

Let $N(\alpha)$ achieve its extremum on $\bar{\alpha}$, and $f(\alpha)$, on
$\tilde{\alpha}$. We have
$$
\begin{aligned}
|N(\tilde{\alpha})-f(\tilde{\alpha})|&<\epsilon,\\
|N(\bar{\alpha})-f(\bar{\alpha})|&<\epsilon.
\end{aligned}
$$
Then
$$
N(\bar{\alpha})<f(\bar{\alpha})+\epsilon<f(\tilde{\alpha})+\epsilon
<N(\tilde{\alpha})+2\epsilon.
$$
But since $N(\alpha)$ is a positive integer, the last inequalities mean that
$$
N(\tilde{\alpha})=N(\bar{\alpha}).
$$
Hence Theorem~\ref{th2} follows.

It can easily be seen that $M(s,n,k)$ increases with $s$. This means that
Theorem~\ref{th2} implies the inequalities
\begin{equation}\label{e67}
\max_{1\le c\le\lfloor ns\rfloor-1}\sum_{i>kc/\lfloor
ns\rfloor}\binom{c}{i}\binom{n-c}{k-i}\le M(s,n,k) \le\max_{1\le c\le\lceil
ns\rceil-1}\sum_{i>kc/\lceil ns\rceil}\binom{c}{i}\binom{n-c}{k-i}.
\end{equation}

\end{document}